\documentclass[a4paper,12pt]{article}
\usepackage{amssymb}
\usepackage{amsmath}
\usepackage{amscd}
\usepackage{amsthm}
\usepackage{graphicx}
\usepackage{enumerate}
\newtheorem{thm}{Theorem}[section]
\newtheorem{lem}[thm]{Lemma}
\newtheorem{rem}[thm]{Remark}

\theoremstyle{definition}
\newtheorem{defi}[thm]{Definition}

\newcommand{\ds}{\displaystyle}

\newcommand{\ep}{\varepsilon}
\newcommand{\s}{\mbox{\boldmath $s$}}
\newcommand{\e}{\mbox{\boldmath $e$}}
\newcommand{\n}{\mbox{\boldmath $n$}}
\renewcommand{\b}{\mbox{\boldmath $b$}}

\renewcommand{\Im}{\operatorname{Im}}
\newcommand{\vect}[1]{\boldsymbol{#1}}
\newcommand{\R}{\boldsymbol{R}}

\newcommand{\rank}{\operatorname{rank}}
\newcommand{\Hess}{\operatorname{Hess}}
\renewcommand{\phi}{\varphi}
\newcommand{\sign}{\operatorname{sgn}}
\newcommand{\inner}[2]{\left\langle{#1},{#2}\right\rangle}

\newcommand{\zv}{\vect{0}}
\renewcommand{\dfrac}{\displaystyle\frac}

\newcommand{\comment}[1]{}
\makeatletter
\renewcommand{\section}{%
  \@startsection{section}%
   {1}%
   {\z@}%
   {-3.5ex \@plus -1ex \@minus -.2ex}%
   {2.3ex \@plus.2ex}%
   {\normalfont\normalsize\bfseries}%
}%
\makeatother
\makeatletter
\renewcommand{\subsection}{%
  \@startsection{subsection}%
   {1}%
   {\z@}%
   {-3.5ex \@plus -1ex \@minus -.2ex}%
   {2.3ex \@plus.2ex}%
   {\normalfont\normalsize\bfseries}%
}%
\makeatother
\begin{document}
\begin{center}
{\bf {\Large Criteria for cuspidal $S_k$ singularities and
their applications}}\\
\bigskip

Kentaro Saji\\
\medskip

\today
\end{center}
 \renewcommand{\thefootnote}{\fnsymbol{footnote}}
 \footnote[0]{Keywords and Phrase. 
 singularities, criteria, Mond singularity, wave front}
\begin{abstract}
 We give useful criteria for $S_1^{\pm}$ singularities
 in the Mond classification table,
 and cuspidal $S_k^{\pm}$ singularities.
 As applications,
 we give a simple proof of a result given by
 Mond and a characterization of cuspidal $S_k^{\pm}$
 singularities for the composition of a cuspidal edge and a
 fold map indicated by Arnol'd for the case $k=1$.
\end{abstract}
\section{Introduction}
Singularities of smooth map-germs have long been studied,
up to the equivalence under coordinate changes in 
both source and target.
There are two separate problems: the classification and the recognition.
The classification is well understood with many good references in
the literature.
Which germ on the classification table is a given germ equivalent to?
Describing simple criteria for this question
 is the recognition and we will do it in this
paper.
In the previous method used for the recognition
a given map-germ is first normalized and then its jet is studied.
The criteria of the recognition without using
normalization are not
only more convenient but also indispensable in
some cases.
We call criteria without normalizing {\it general criteria\/}
for a while.
In fact, in the case of wave front surfaces in 3-space,
general criteria for the cuspidal edge and the swallowtail
were
given in \cite{krsuy}
where we studied the local and global behavior of flat fronts in
hyperbolic 3-space using them.
Moreover 
the singular curvature
on the cuspidal edge was introduced
and its properties were investigated in \cite{SUYfr}.
Furthermore, a general criterion for the cuspidal cross cap
was
given in \cite{FSUY},
where we studied maximal surfaces and 
constant mean curvature one surfaces in the Lorentz-Minkowski 3-space
and described a certain duality between the swallowtails
and the cuspidal cross caps.
The cuspidal cross cap is also called the {\it cuspidal\/} $S_{0}$
{\it singularity}.
In \cite{ist}, general criteria for 
the cuspidal lips and the cuspidal beaks
were given and the horo-flat surfaces in hyperbolic space
were investigated.
Recently, several applications of
these criteria were considered in various situations
\cite{ishi-machi,isman,ist,circular,kruy,sch}.
Criteria for higher dimensional $A$-type singularities of wave fronts
and their applications
were considered in \cite{SUYcam}.

\begin{figure}[htbp]
\centering
\begin{tabular}{ccc}
   \includegraphics[width=.3\linewidth]{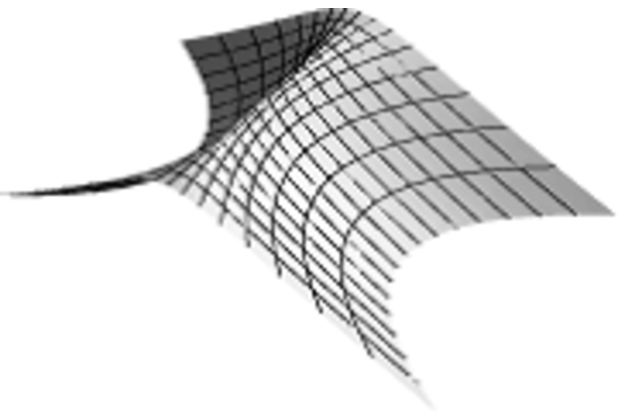}&
   \hspace*{10mm}&
   \includegraphics[width=.2\linewidth]{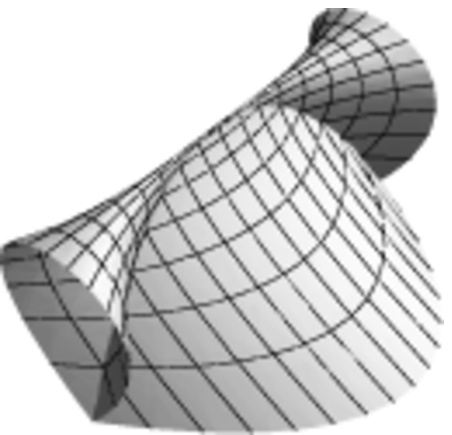}
   \end{tabular}
 \caption{The Chen Matumoto Mond $+$ singularity (left) and
  $-$ singularity (right).}
 \label{fig:cmm}
 \end{figure}

In this paper,
we first give general criteria 
for the {\it Chen Matumoto Mond\/ $\pm$
singularities\/} $S_{1}^\pm$ which are
map-germs defined by
\begin{equation}
\label{eq:cmm}
 S^{\pm}_1:(x,y)\mapsto\big(x,y^2,y(x^2\pm y^2)\big)
\end{equation}
at the origin (See Figure \ref{fig:cmm}).

\begin{figure}[htbp]
 \centering
 \begin{tabular}{ccccc}
  \includegraphics[width=.2\linewidth]{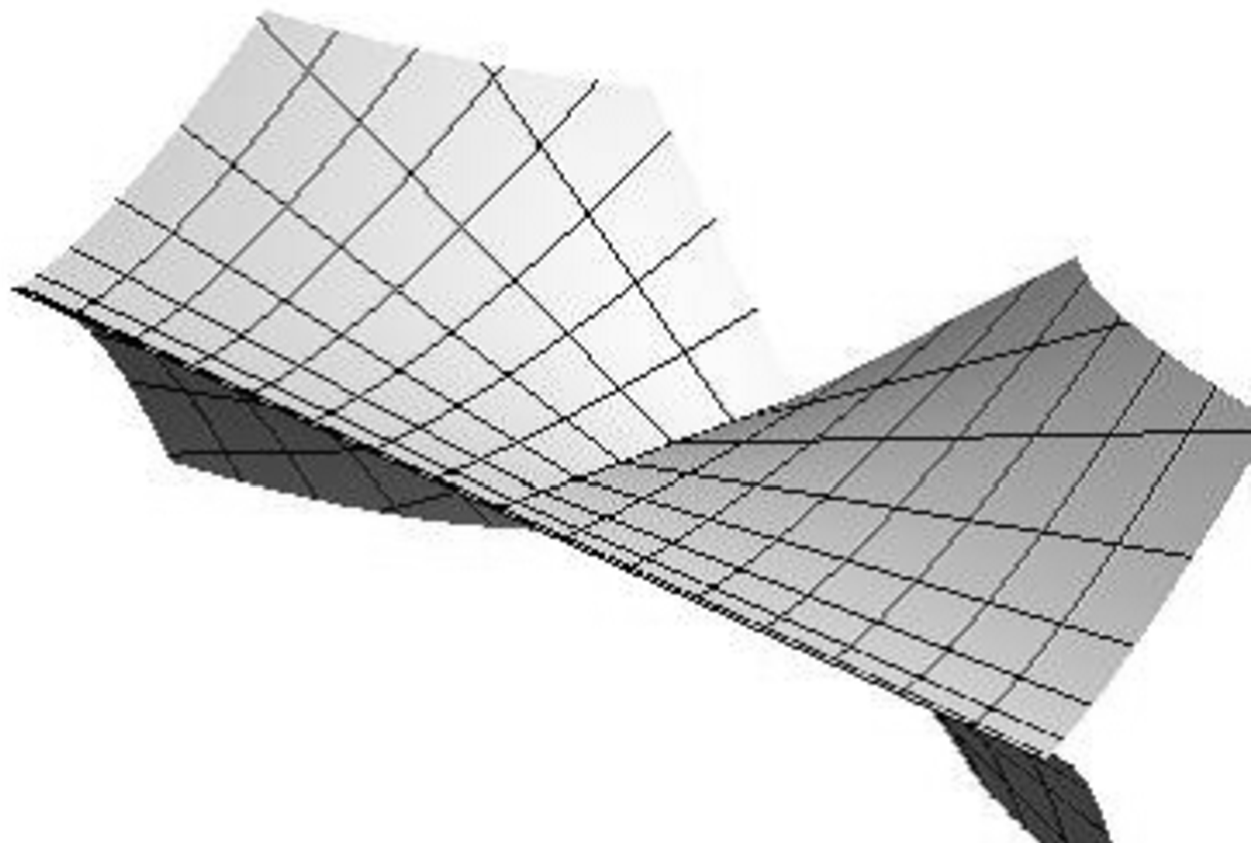}&
  \hspace*{10mm}&
  \includegraphics[width=.17\linewidth]{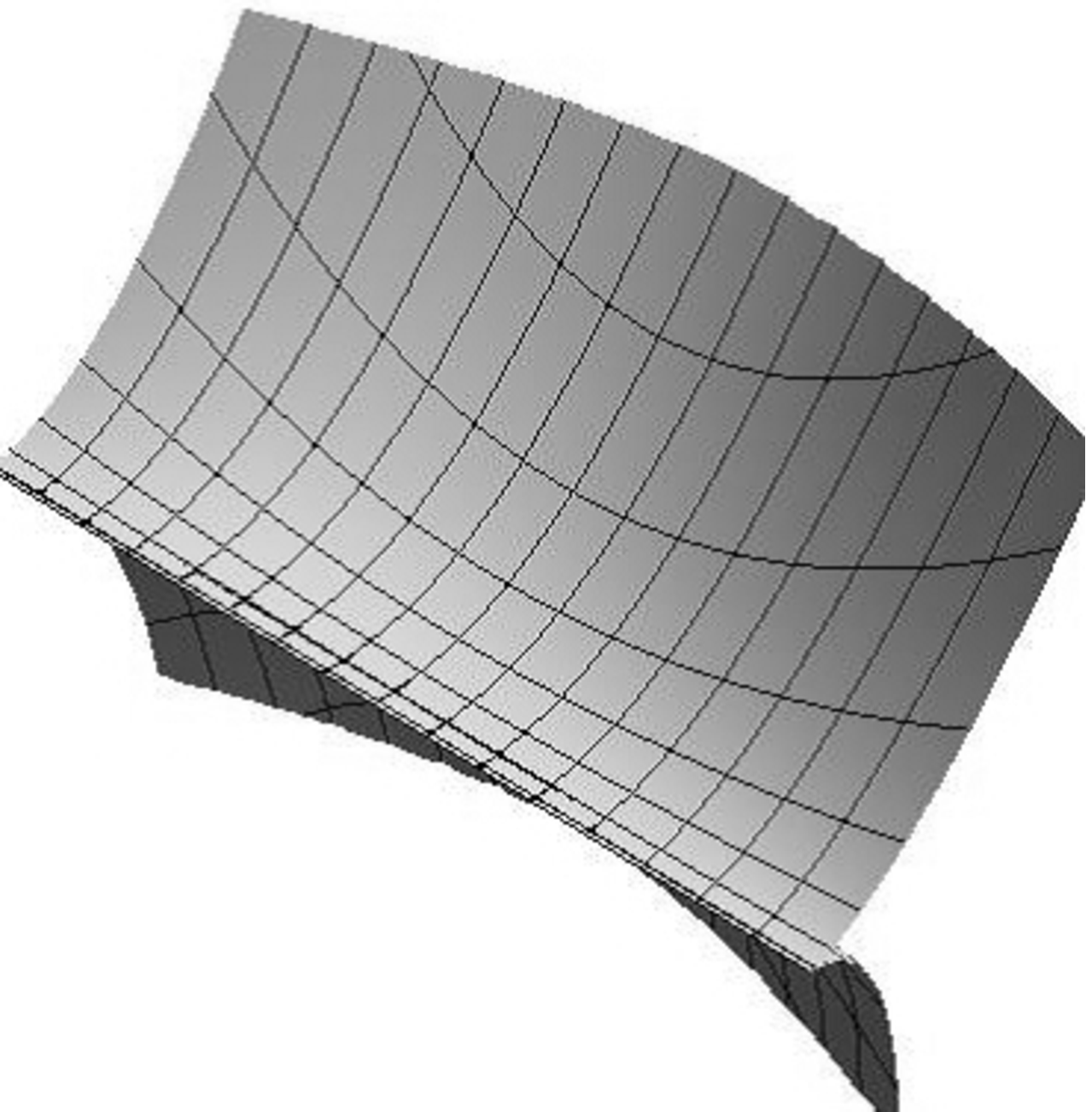}
  \hspace*{10mm}&
  \includegraphics[width=.2\linewidth]{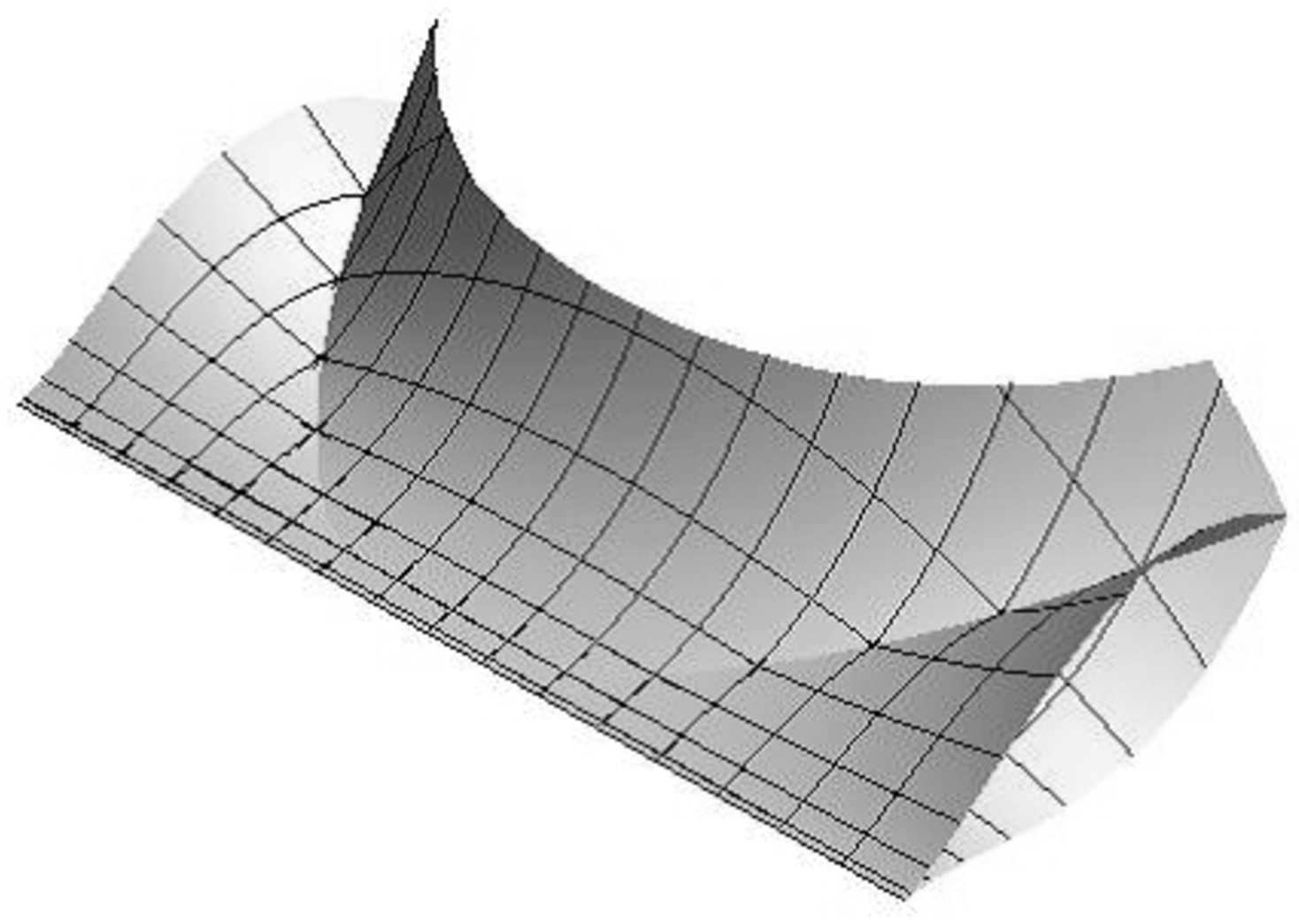}
 \end{tabular}
 \caption{$cS_{0}$ singularity (left), 
  $cS_{1}^+$
 singularity (center) and
 $cS_{1}^-$ singularity (right).}
 \label{fig:csk}
\end{figure}

X. Y. Chen and T. Matumoto showed these singularities and their
suspensions are the generic singularities
of one-parameter families of $n$-dimensional
manifolds in $\R^{2n+1}$ (\cite{chenmatu}).
In \cite{Mond2}, D. Mond classified simple singularities
$\R^2\to\R^3$ with respect to the ${\cal A}$-equivalence,
giving a criterion for
map-germs of the forms
$(x,y)\mapsto(x,y^2,f(x,y))$ \cite[Theorem 4.1.1]{Mond2}.
The Chen Matumoto Mond $\pm$ singularities
appear as $S_{1}^\pm$ singularities
in his classification
table \cite{Mond2}.
In this paper, we also give criteria for the {\it cuspidal\/ $S_k^\pm$ 
singularities},
which are map-germs defined by
\begin{equation*}
 cS^{\pm}_k:(x,y)\mapsto\big(x,y^2,y^3(x^{k+1}\pm y^2)\big),
 \qquad (k=0,1,\ldots)
\end{equation*}
at the origin (See Figure \ref{fig:csk}).
These are kinds of ``cusped'' $S_k^{\pm}$ singularities.
If $k$ is even,
the cuspidal $S_k^+$ singularity
 and the cuspidal $S_k^-$ singularity are ${\cal A}$-equivalent.
If $k=0$, this is the {\it cuspidal cross cap}.
We state criteria for the cuspidal $S_k^\pm$ singularities
as a generalization of the criterion for the cuspidal cross cap
given in \cite{FSUY}.
It is known that the
 cuspidal $S_k^\pm$ singularities appear as singularities
of frontal surfaces (for the definition of frontal surfaces, see \S
3).
As applications,
we give a simple proof of the properties on
singularities of tangent developable surfaces
given by Mond \cite{Mond1} and
an interpretation of the degree of contactness
about V. I. Arnol'd's observation \cite{arnld-sing-caus}
in \S \ref{sec:app}.
All maps considered here are of class $C^\infty$.

The author would like to express his sincere gratitude
to Professor Takashi Nishimura for his invaluable advice
and comments.
In particular, he had suggested a proof of Theorem \ref{thm:csk}
using Thom's splitting lemma simpler than the
preliminary version.
This work has partially done during the author's stay at
the University of Sydney, as a participant of the
JSPS joint research program with Australia 2008.
The author would like to thank Professors
Laurentiu Paunescu and Satoshi Koike for fruitful discussions and
kind hospitality.
The author is grateful to the referee 
for careful reading and helpful comments.

\section{Criteria for the Chen Matumoto Mond $\pm$ singularities}
In this section, we show criteria for the Chen Matumoto Mond
singularities
of surfaces.
If a map-germ $f:(\R^2,\zv)\to(\R^3,\zv)$ satisfies
$\rank df_{\zv}=1$, the singular point $\zv$ is called {\it corank one}.
If $f:(\R^2,\zv)\to(\R^3,\zv)$ has a corank one singular point at
$\zv$,
then there exist vector fields $(\xi,\eta)$
near the origin such that $df_{\zv}(\eta_{\zv})=\zv$ and
$\xi_{\zv},\eta_{\zv}\in T_{\zv}\R^2$ are linearly independent. 
We define a function $\phi:(\R^2,\zv)\to\R$ by
\begin{equation}
 \label{eq:cr}  
  \phi=\det(\xi f,\ \eta f,\ \eta\eta f),
\end{equation}
where $\zeta g:(\R^2,\zv)\to(\R^3,\zv)$ is the directional derivative of
a
vector valued
function $g$ by a vector field $\zeta$.
We call $\eta_{\zv}$ the {\it null direction\/} (cf. \cite{krsuy}).
\begin{defi}
  Two map-germs $f_i:(\R^2,\zv)\to(\R^3,\zv)$ $(i=1,2)$ are
  {\it ${\cal A}$-equivalent\/} if there exist
  germs of $C^\infty$-diffeomorphisms
  $d_s:(\R^2,\zv)\to(\R^2,\zv)$
  and $d_t:(\R^3,\zv)\to(\R^3,\zv)$ such that
  $d_t\circ f_1=f_2\circ d_s$
  holds.
\end{defi}

\begin{thm}
 \label{thm:cmm}
 Let\/ $f:(\R^2,\zv)\to(\R^3,\zv)$ be a map-germ and\/
 $\zv$ a corank one singular point.
 Then\/ $f$ at\/ $\zv$ is\/ ${\cal A}$-equivalent
 to the Chen Matumoto Mond\/ $-$
 singularity if and only if\/
 $\phi$ has a critical point at\/ $\zv$ and\/
$\det\Hess\phi(\zv)>0$.
 On the other hand\/ $f$ at\/ $\zv$ is\/ ${\cal A}$-equivalent to the
 Chen Matumoto Mond\/ $+$
 singularity
 if and only if\/
 $\phi$ has a critical point at\/ $\zv$,
 $\det\Hess\phi(\zv)<0$ and\/  
 two vectors\/ $\xi f(\zv)$ and\/ 
 $\eta\eta f(\zv)$ are linearly
 independent.
\end{thm}
\begin{rem}
 \begin{itemize}
  \item  The additional condition in the case\/ $\det\Hess\phi<0$ cannot be
         removed.
         For example, $(x,xy+y^3,xy+2y^3)$ satisfies\/
         $\det\Hess\phi(\zv)<0$
         but\/ $\xi f(\zv)$ and\/ $\eta\eta f(\zv)$ are linearly dependent.
         It is known that this map-germ is
         not\/ ${\cal A}$-equivalent to the Chen Matumoto Mond
         singularity.
         If\/ $\det\Hess\phi>0$, then\/
         $\xi f(\zv)$ and\/ $\eta\eta f(\zv)$ are automatically
         linearly independent.
  \item Using the above function\/ $\phi$, we can write the recognition
        criterion for Whitney umbrella by\/
        $\xi \phi\ne0$, that is, $d \phi\ne\zv$.
  \item Since\/ $\eta f(\zv)=\zv$, Theorem\/ \ref{thm:cmm} implies
        that the Chen Matumoto Mond singularity is three determined,
        namely, if the\/ $3$-jet of a map-germ\/
        $f:(\R^2,\zv)\to(\R^3,\zv)$
        is\/ ${\cal A}$-equivalent to
        the Chen Matumoto Mond singularity, then\/
        $f$ is\/ ${\cal A}$-equivalent to
        the Chen Matumoto Mond singularity.
 \end{itemize}
\end{rem}
To prove Theorem \ref{thm:cmm},
the following lemmas play the key role.
\begin{lem}
 \label{lem:source}
 The conditions in Theorem\/ {\rm \ref{thm:cmm}} are
 independent of the choice of vector fields\/ $(\xi,\eta)$.
\end{lem}
\begin{lem}
 \label{lem:target}
 The conditions in Theorem\/ {\rm \ref{thm:cmm}} are
 independent of the choice of coordinates on the target.
\end{lem}
\begin{proof}[Proof of Lemma {\rm \ref{lem:source}}.]
Let us put
 $$
 \left\{
 \begin{array}{rcl}
  \overline{\xi}&=&a_{11}\xi+a_{12}\eta\\
  \overline{\eta}&=&a_{21}\xi+a_{22}\eta
 \end{array}
 \right.
 ,\quad
 \big((a_{ij}):\R^2\to GL(2,\R),\ a_{21}(\zv)=0\big),
 $$
 and
 $$
 \overline{\phi}=\det(\overline{\xi}f,\
 \overline{\eta}f,\ \overline{\eta}\,\overline{\eta}f).
 $$
 Then by a straight calculation,
 we have
 $$
 \begin{array}{rcl}
  \overline{\xi}f&= &a_{11}\xi f+a_{12}\eta f, \\
  \overline{\eta}f&= &a_{21}\xi f+a_{22}\eta f \text{\ \ and\ } \\
  \overline{\eta}\,\overline{\eta}f&= &
   *\xi f+
   *\eta f+a_{21}(a_{21}\xi\xi f+a_{22}\xi\eta f+a_{22}\eta\xi f)
   +a_{22}^2\eta\eta f.
 \end{array}
 $$
 Thus it follows that the linear independence of two
 vectors
 $\xi f(\zv)$ and $\eta\eta f(\zv)$ does not depend on the choice of
 vector fields.
 Hence, we have
 $$
 \overline{\phi}=(a_{11}a_{22}-a_{12}a_{21})
 \Big(a_{21}\det(\xi f,\ \eta f,\ a_{21}\xi\xi f+
 a_{22}\xi\eta f+a_{22}\eta\xi f)
 +a_{22}^2\det(\xi f,\ \eta f,\ \eta\eta f)\Big).
 $$
 Now it is sufficient to prove that 
 $$
 \begin{array}{l}
 \xi m(\zv)=\eta m(\zv)=0\ \text{and}\ \Hess m(\zv)=O,\ \\
 \hspace{40mm}\text{where}\ 
 m:=
 a_{21}\det(\xi f,\ \eta f,\ a_{21}\xi\xi f+
 a_{22}\xi\eta f+a_{22}\eta\xi f).
 \end{array}
 $$
 Since $m$ contains the terms $a_{21}$ and $\eta f$,
 which vanish at the origin,
 it holds that $\xi m(\zv)=\eta m(\zv)=0$.
 Next, we have assumed that $\phi$ has a critical point
 at $\zv$, namely,
 \begin{equation}
  \label{eq:criti}
   \xi\phi(\zv)=\det(\xi f,\ \xi\eta f,\ \eta\eta f)(\zv)=0.
 \end{equation}
 Since $\xi\eta-\eta\xi$ is a vector field, and $\xi$ and $\eta$ are
 linearly independent,
 $\xi\eta-\eta\xi$ is a linear combination of $\xi$ and $\eta$
 at each point.
 Moreover,
 $\xi\eta f(\zv)-\eta\xi f(\zv)$ is parallel to $\xi f(\zv)$,
 since $\eta$ is the null vector at $\zv$.
 Thus we see that
 $-(\xi\eta+\eta\xi)f=-2\xi\eta f+(\xi\eta-\eta\xi)f$
 is a linear combination of $\xi\eta f$ and $\xi f$ at $\zv$.
 Thus,
 $\xi\xi m(\zv)=2\xi a_{21}(\zv)
 \det(
 \xi f,\ \xi\eta f,\ a_{21}\xi\xi f+
 a_{22}\xi\eta f+a_{22}\eta\xi f)(\zv)=0$ holds,
 since $a_{21}(\zv)=0$.
 By the same reason and (\ref{eq:criti}), we also have
 $$
 \begin{array}{l}
  \xi\eta m(\zv)=\xi a_{21}(\zv)
  \det(
  \xi f,\ \eta\eta f,\ a_{21}\xi\xi f+
  a_{22}\xi\eta f+a_{22}\eta\xi f)(\zv)
  \\ \hspace{30mm}+
  \eta a_{21}(\zv)
  \det(
  \xi f,\ \xi\eta f,\ a_{21}\xi\xi f+
  a_{22}\xi\eta f+a_{22}\eta\xi f)(\zv)
  =0.
 \end{array}
 $$
Furthermore,
 $\eta\eta\,m(\zv)=2\eta a_{21}(\zv)
 \det(
 \xi f,\ \eta\eta f,\ a_{21}\xi\xi f+
 a_{22}\xi\eta f+a_{22}\eta\xi f)(\zv)
 =0$ as well.
 Hence $\Hess m(\zv)=O$ holds.
\end{proof}
\begin{proof}[Proof of Lemma {\rm \ref{lem:target}}.]
 Take a $C^\infty$-diffeomorphism
 $\Phi(X)=(\Phi_1(X),\Phi_2(X),\Phi_3(X)):\R^3\to\R^3$, where
 $X=(X_1,X_2,X_3)$.
 Put $f=(f_1,f_2,f_3)$ and
 $$
 \tilde{\phi}=\det\Big(\xi(\Phi\circ f),\ \eta(\Phi\circ f),\
 \eta\eta(\Phi\circ f)\Big).
 $$
 Then 
 the first component of
 the vector 
 $\eta\eta(\Phi\circ f)=\eta\big(d\Phi(\eta f)\big)$
 is calculated as
 \begin{equation}
  \label{eq:etaetaphif}
 \begin{array}{l}
  \eta
 \left(
   \ds\sum_{i=1}^3
   \dfrac{\partial \Phi_1}{\partial X_i}\eta f_i
 \right)
 =
 \ds\sum_{i=1}^3\left(\left(
   \ds\sum_{j=1}^3
   \dfrac{\partial^2 \Phi_1}{\partial X_i\partial X_j}\eta f_j
 \right)
 \eta f_i
 +\dfrac{\partial \Phi_1}{\partial X_i}\eta\eta f_i\right)\\[6mm]
 \hspace{80mm}=
 \Hess\Phi_1(\eta f,\eta f)+d\Phi_{1}(\eta\eta f)
 .
 \end{array}
 \end{equation}
 Hence the linear independence of
 $\xi f(\zv)$ and $\eta\eta f(\zv)$ 
 does not depend on the choice of the coordinates of the target.
 By (\ref{eq:etaetaphif}) again, it holds that 
 $$
 \begin{array}{rcl}
  \tilde{\phi}&=&
     \det\left(d\Phi(\xi f),\ d\Phi(\eta f),\ 
\left(
\begin{array}{c}
 \Hess\Phi_1(\eta f,\eta f)\\
 \Hess\Phi_2(\eta f,\eta f)\\
 \Hess\Phi_3(\eta f,\eta f)\\
\end{array}   {}\right)\right)
   +(\det d\Phi)\,\phi.
 \end{array}
 $$
  Thus by the same argument as above, it is sufficient to prove that
 $\Hess M(\zv)=O$,
 where
 $$M:=\det\Big(d\Phi(\xi f),\ d\Phi(\eta f),\
 \Hess\Phi_i(\eta f,\eta f)_{i=1,2,3}\Big).$$
 Since $\eta f$ vanishes at the origin, $\Hess M(\zv)=O$ holds.
\end{proof}
Using these Lemmas, we prove Theorem \ref{thm:cmm}.
\begin{proof}[Proof of Theorem\/ {\rm \ref{thm:cmm}}.]
 The necessity of the conditions
 is immediate from the calculation for the formula (\ref{eq:cmm})
 and Lemmas \ref{lem:source} and \ref{lem:target}.
 We prove that the conditions are sufficient.
 Let us assume the conditions in Theorem \ref{thm:cmm}.
 By Lemmas \ref{lem:source} and \ref{lem:target},
 we can change vector fields $(\xi,\eta)$ and
 coordinates on the target.
 Moreover, since the conditions do not depend on the coordinates on
 the source, we may change coordinates on the source.
 Since $f$ is corank one at $\zv$,
 by the implicit function theorem,
 $f$ is ${\cal A}$-equivalent to the map-germ defined by
 $(x,y)\mapsto(x,f_2(x,y),f_3(x,y))$ at the origin.
 By the target coordinate change,
 $f$ is ${\cal A}$-equivalent to the map-germ $(x,yg(x,y),yh(x,y))$.
 Since $f$ has a singularity at the origin,
 there is no constant term in $g$ and $h$.
 Moreover, we have the following lemma.
 \begin{lem}
   \label{lem:pr1}
   At the origin,
   $g_y$ or\/ $h_y$ does not vanish, where\/ $g_y=\partial g/\partial y$
   and\/ $h_y=\partial h/\partial y$.
\end{lem}
 \begin{proof}
  Since $(\partial/\partial y)\big(x,yg(x,y),yh(x,y)\big)(\zv)=\zv$ holds,
  we may choose $\xi=\partial/\partial x,\ \eta=\partial/\partial y$.
  Then it holds that 
  $$
  \phi=\det\left(
  \begin{array}{ccc}
   1&0&0\\
   *&g+yg_y&2g_y+yg_{yy}\\
   *&h+yh_y&2h_y+yh_{yy}
  \end{array}
  \right).
  $$
  Since $g(\zv)=h(\zv)=0$, we have
  $\phi_{yy}(\zv)=6(g_yh_{yy}-g_{yy}h_y)(\zv)$.

  In the case of $\det\Hess\phi(\zv)>0$, if we assume that 
  $g_y(\zv)=h_y(\zv)=0$, then $\phi_{yy}(\zv)=0$
  holds and hence $\det\Hess\phi(\zv)
  =-(\phi_{xy})^2(\zv)\leq0$ at the origin,
  which contradicts the assumption.
  Therefore, in this case
  $g_y(\zv)\ne0$ or $h_y(\zv)\ne0$ holds.

  On the other hand, in the case of $\det\Hess\phi<0$,
  we have the additional condition which implies that
  $\eta\eta f(\zv)\ne\zv$.
  Thus we have $g_y(\zv)\ne0$ or $h_y(\zv)\ne0$.
 \end{proof}
 Let us continue the proof of Theorem~\ref{thm:cmm}.
 By Lemma \ref{lem:pr1}, we may assume
 $g_y(\zv)\ne0$.
 Since by the implicit function theorem,
 the set $\{g(x,y)=0\}$ is a regular curve,
 we take a new coordinate system 
 $(x,\tilde y)$ satisfying
 $g(x,\tilde y)=0$ on $\tilde y=0$.
 Then we may assume that
 $f(x,\tilde y)=
 \big(x,\tilde y^2g(x,\tilde y), \tilde yh(x,\tilde y)\big)$,
 $g(0,0)\ne0$.
 Furthermore, considering a coordinate change
 $(\overline{x},\overline{y})
 =
 \big(x,\tilde y\sqrt{|g(x,\tilde y)|}\big)$
 and rewriting $(\overline{x},\overline{y})$ by
 $(x,y)$, we may assume 
 $
 f(x,y)=(x,y^2,y\overline{h}(x,y)).$
 Now we set
 $$
 \overline{h}_1(x,y)=\frac{\overline{h}(x,y)+\overline{h}(x,-y)}{2}
 ,\qquad
 \overline{h}_2(x,y)=\frac{\overline{h}(x,y)-\overline{h}(x,-y)}{2}.
 $$
 Then $\overline{h}(x,y)=\overline{h}_1(x,y)+\overline{h}_2(x,y)$
 holds and $\overline{h}_1(x,y)$ (resp. $\overline{h}_2(x,y)$) is
 an even (resp. odd) function with respect to $y$.
 Then there exist
 functions $\tilde{h}_1(x,y)$ and $\tilde{h}_2(x,y)$
 such that
 $$
 \overline{h}_1(x,y)=\tilde{h}_1(x,y^2),\qquad
 \overline{h}_2(x,y)=y\tilde{h}_2(x,y^2).
 $$
 Then we have
 $f(x,y)=\big(x,y^2,y\tilde{h}_1(x,y^2)+y^2\tilde{h}_2(x,y^2)\big)$.
 Considering a coordinate change
 $\tilde\Theta(X,Y,Z)=(X,Y,Z-Y\tilde{h}_2(X,Y))$
 and replacing $f$ by $\tilde\Theta\circ f$, we may set
 $$f(x,y)=\big(x,y^2,y\tilde{h}_1(x,y^2)\big).$$
 Since the function $\phi$ defined by \eqref{eq:cr}
 for this map has the form
 $
 -2\tilde{h}_1(x,y^2)+y\,*$, 
 it holds that
 $(\partial/\partial x)\tilde{h}_1(\zv)$ $=0$.
 Here $*$ means a function.
 Thus
 there exists a function $\tilde{f}$ such that
 $$
 f(x,y)=\bigg(x,y^2,y\Big[\alpha x^2+\beta
 y^2\big(1+\tilde{f}(x,y^2)\big)\Big]\bigg),\quad
 \tilde{f}(0,0)=0.
 $$
 Note that the function $\phi$ for this map has the form
 $$
 -2\alpha x^2+6\beta y^2+\text{(higher\ order\ term)}.
 $$
 Considering a diffeomorphism $\theta$ defined by
 \begin{equation*}
   (u,v)=\theta(x,y)=(x,\theta_2(x,y))=
 \left(x,y\sqrt{1+\tilde{f}(x,y^2)}\right)
 \end{equation*}
 and the inverse map $\theta^{-1}(u,v)=(u,\vartheta_2(u,v))$,
 we have the following lemma.
 \begin{lem}
  \label{lem:thetainv}
  There exists a function\/ $f_2$ satisfying\/
  $f_2(\zv)\ne0$ and\/
  $\vartheta_2(u,v)=vf_2(u,v^2)$.
 \end{lem}
 \begin{proof}
  Substituting $v=0$ in the identity
  \begin{equation*}
   \label{eq:theta}
    \theta\circ\theta^{-1}(u,v)
    =
    \left(u,\vartheta_{2}(u,v)
     \sqrt{1+\tilde{f}(u,\vartheta_{2}(u,v)^2)}\right)
    =(u,v),
  \end{equation*}
  we have $\vartheta_{2}(u,0)=0$.
  Next, we take $(x,y)=\theta^{-1}(u,-v)$.
  Then we have $v=-\theta_2(x,y)=\theta_2(x,-y)$.
  Since $(u,v)=(x, \theta_2(x,-y))$, it holds that
  $\theta^{-1}(u,v)=(x,-y)$.
  Thus $\vartheta_2(u,-v)=-\vartheta_{2}(u,v)$.
  Hence $\vartheta_2$ satisfies that
  $\vartheta_2(u,0)=0$ and $\vartheta_2(u,-v)=-\vartheta_2(u,v)$
  for any $(u,v)$ near $\zv$.
  Then by the same argument as the construction of
  $\bar{h}_2(x,y)$ in the proof of Theorem \ref{thm:cmm},
  we have the lemma.
 \end{proof}
 By Lemma \ref{lem:thetainv}, 
 the composition $f\circ \theta^{-1}$ has the expression
 $$
 \big(u,v^2f_2(u,v^2)^2,
 vf_2(u,v^2)(\alpha u^2+\beta v^2)\big).
 $$
 Considering a diffeomorphism
 $\Theta(X,Y,Z)=(X,Yf_2(X,Y)^2,Zf_2(X,Y))$,
 we see that $\Theta^{-1}\circ f\circ \theta^{-1}$
 has the form
 $(u,v^2,v(\alpha u^2+\beta v^2))$.
 This is ${\cal A}$-equivalent to the desired map-germ
 because we see that $$-\sign(\alpha\beta)=\sign(\det\Hess\phi(\zv)).$$
\end{proof}
\section{Criteria for cuspidal $S_k^\pm$ singularities of frontals}
In this section, we shall introduce
the notion of frontal surfaces
and give criteria for the cuspidal $S_k^\pm$ singularities 
of frontals.
\subsection{Preliminaries on the frontals}
The projective cotangent bundle $PT^*\R^3$ of  $\R^3$
has the canonical contact structure
and can be identified with the projective tangent bundle $PT \R^3$.
A smooth map-germ $f:(\R^2,\zv) \to (\R^3,\zv)$ is called a {\it frontal\/}
if there exists a never-vanishing vector field $\nu$ of $\R^3$ along $f$
such that $L:=(f,[\nu]):(\R^2,\zv)\to (\R^3\times
P^2,L(\zv))=(PT\R^3,L(\zv))$
is an isotropic map,
that is, the pull-back of the canonical contact form of $PT\R^3$
vanishes on $\R^2$,
where $P^2$ means the projective space
and $[\nu]$ means the projective class of $\nu$.
This condition is equivalent to the following orthogonality condition:
\begin{equation*}
   \inner{df(X_p)}{\nu(p)}=0 \qquad
    (^\forall p\in\R^2,\quad{^\forall} X_p\in T_p\R^2),
\end{equation*}
where $\inner{~}{~}$ is the canonical inner product on $\R^3$.
The vector field  $\nu$ is called the {\it normal vector\/} of
the frontal $f$.
The plane perpendicular to $\nu(p)$ is called the {\it limiting tangent
plane\/}
at $p$.
A frontal $f$ is called a {\it front\/} if
$L=(f,[\nu])$ is an immersion 
(cf. \cite{arnld-sing-caus} see also \cite{krsuy}).
A function
\begin{equation*}
   \lambda(u,v):=\det(f_u,f_v,\nu)
\end{equation*}
is called the {\it signed area density function},
where $(u,v)$ is the coordinate system on $\R^2$.

Let $\zv$ be a singular point of a frontal
$f:(\R^2,\zv)\to(\R^3,\zv)$.
Then the set of singular points $S(f)$ of $f$ coincides with
the zeros of $\lambda$ near $\zv$.
If $d\lambda(\zv)\ne0$, then $\zv$ is called a {\it non-degenerate 
singular point}.
Assume now that
$\zv$ is a non-degenerate singular point.
Then there exists a regular curve
$\gamma(t):((-\ep,\ep),0)\to(\R^2,\zv)$ $(\ep>0)$
such that the image of $\gamma$ is $S(f)$.
Also
dimension of the kernel $\ker(df_{\gamma(t)})$ is equal to one
and there is
a never-vanishing vector $\eta(t)$ such that
$\eta(t)$ spans $\ker(df_{\gamma(t)})$.
We call $\eta$ the {\it null vector field}.
We define a function $\psi$ on $S(f)$ by
\begin{equation}
 \label{eq:ccr}
  \psi(t)=\det\left(\frac{df\circ\gamma}{dt}(t),
               \nu\circ\gamma(t),
               d\nu_{\gamma(t)}\big(\eta(t)\big)\right)\quad
  \text{for}\quad t\in(-\varepsilon,\varepsilon).
\end{equation}
This function is originally defined in \cite{FSUY}.
The signed area density function, the non-degeneracy and
the null vector field are introduced in \cite{krsuy}.
\subsection{Criterion for the $(2,5)$-cusp}
If we substitute $u=0$ in the normal form of the 
cuspidal $S_k^{\epsilon}$ singularity
$\big(u,v^2,v^3(u^{k+1}+\epsilon v^2)\big)$, $\epsilon=\pm1$,
it reduces to a
$(2,5)$-cusp curve through $\zv$.
In this subsection, we state a criterion for the $(2,5)$-cusp,
namely, the map-germ
given by $t\mapsto(t^2,t^5,0)$ at $t=0$.
\begin{lem}\label{lem:25}
 Let\/ $c(t):I\to\R^3$ be a curve and\/
 $0\in I$.
 
 {\rm (i)} Assume that\/ $c$ satisfies\/ $c'_0=\zv$, $c''_0\ne\zv$
 ,$c'''_0=c^{(4)}_0=\zv$, and\/ $c''_0$ and\/ $c^{(5)}_0$ are linearly
 independent,
 then\/ $c$ at\/ $t=0$ is\/ ${\cal A}$-equivalent to the\/
 $(2,5)$-cusp.

 {\rm (ii)} Assume that\/ $c$ satisfies\/ 
 $c'_0=\zv$, $c''_0\ne0$ and
 two vectors\/ $c''_0$ and\/ $c'''_0$ are linearly dependent,
 that is, there exists\/ $\ell\in\R$ such that\/
 $c'''_0=\ell c''_0$.
 If 
 two vectors\/ $c''_0$ and\/ $3c^{(5)}_0-10\ell c^{(4)}_0$ are linearly
 independent in addition,
 then\/ $c$ at\/ $t=0$ is\/ ${\cal A}$-equivalent to the\/
 $(2,5)$-cusp.
\end{lem}
 Here,
 $c'=c^{(1)}=dc/dt$, $c^{(j)}=dc^{(j-1)}/dt$,
 and $c^{(j)}_0=c^{(j)}(0)$ $(j=1,\ldots,5)$.
\begin{proof}
 One can prove {\rm (i)} by a fundamental argument.
 So, we omit its proof. We shall prove (ii).
 Suppose that $c$ satisfies the assumptions of (ii)
 except the last condition.
 Then $c$ is written as
 $$
 \begin{array}[tb]{l}
  c(t)=
   \big(
   a_2t^2+a_3t^3+a_4t^4+a_5t^5+o(t^5),
   ka_2t^2+ka_3t^3+b_4t^4+b_5t^5+o(t^5),\\
  \hspace{15mm}ka_2t^2+ka_3t^3+c_4t^4+c_5t^5+o(t^5)
   \big),\ 
   a_2,a_3,a_4,a_5,b_4,b_5,c_4,c_5,k\in\R \text{\ and\ }
   a_2\ne0,
 \end{array}
 $$
 where, $o(t^5)$ is a Landau notation.
 Considering a coordinate change on the target
 $(X,Y,Z)$ $\mapsto(X,Y-kX,Z-kX)$, we see that $c$ is ${\cal A}$-equivalent
 to
 $$
 \begin{array}{l}
 \big(
 a_2t^2+a_3t^3+a_4t^4+a_5t^5+o(t^5),
 (b_4-ka_4)t^4+(b_5-ka_5)t^5+o(t^5),\\
 \hspace{70mm}(c_4-ka_4)t^4+(c_5-ka_5)t^5+o(t^5)
 \big).
 \end{array}
 $$
 Next, considering a parameter change $t\mapsto t-(a_3/2a_2)t^2$, we get
 $$
 \begin{array}{l}
  \bigg(
   a_2t^2+
   \left(-\ds\frac{5a_3^2}{4a_2}+a_4\right)t^4
   +\left(\ds\frac{3a_3^3}{4a_2^2}-
   \ds\frac{2a_3a_4}{a_2}+a_5\right)t^5+o(t^5),\\
  \hspace{20mm}(b_4-ka_4)t^4+
   \left(-\ds\frac{2a_3}{a_2}(b_4-ka_4)+b_5-ka_5\right)t^5+o(t^5),\\
  \hspace{40mm}(c_4-ka_4)t^4+
   \left(-\ds\frac{2a_3}{a_2}(c_4-ka_4)+c_5-ka_5\right)t^5+o(t^5)
   \bigg).
 \end{array}
 $$
 Lastly, considering a coordinate change 
 $(X,Y,Z)\mapsto
 (X-(-5a_3^2+4a_2a_4)X^2/(4a_2^2),
 Y-(b_4-ka_4)X^2/a_2^2,Z-(c_4-ka_4)X^2/a_2^2),$
 we get
 $$
 \begin{array}[tb]{l}
  \bigg(
   a_2t^2
   +\left(\ds\frac{3a_3^3}{4a_2^2}-\ds\frac{2a_3a_4}{a_2}+a_5\right)t^5+o(t^5),
   \left(-\ds\frac{2a_3}{a_2}(b_4-ka_4)+b_5-ka_5\right)t^5+o(t^5),\\
  \hspace{80mm}
   \left(-\ds\frac{2a_3}{a_2}(c_4-ka_4)+c_5-ka_5\right)t^5+o(t^5)
   \bigg).
 \end{array}
 $$
 By a direct calculation we see the last condition of {\rm (ii)}
 is equivalent to the condition
 $$
 a_2(b_5-ka_5)-2a_3(b_4-ka_4)\ne0
 \text{\ or\ }
 a_2(c_5-ka_5)-2a_3(c_4-ka_4)\ne0\quad \text{at}\quad t=0.
 $$
 This is also the assumption of {\rm (i)} for
 the curve with respect to the last coordinate
 change and we complete the proof.
\end{proof}
\comment{Remark that
by this lemma, the linearly
independence of two vectors $c''_0$ and $3c^{(5)}_0-10\ell c^{(4)}_0$
depends neither on the choice of parameter $t$ nor the
coordinates on the target space for $\ell\in\R$
satisfying $c'''_0=\ell c''_0$.}

\subsection{Criteria for cuspidal $S_k^\pm$ singularities}
Criteria for the cuspidal $S_k^\pm$ singularities are
stated as follows:
\begin{thm}\label{thm:csk}
 Let\/ $f:(\R^2,\zv)\to(\R^3,\zv)$ be a frontal and
 fix a representative\/ $\nu$ of the normal vector of\/ $f$.
 The map-germ\/ $f$ at\/ $\zv$ is\/ ${\cal A}$-equivalent
 to the cuspidal\/ $S^\pm_{k-1}$ 
 singularity\/ $(k\ge2)$
 if and only if the following\/ {\rm (a)-(d)} hold\/{\rm :} 
 \begin{enumerate}[{\rm (a)}]
  \item\label{it:ndg}
       $\zv$ is a non-degenerate singular point
       and the null vector is transverse to\/ $S(f)$ at\/ $\zv$.
  \item\label{it:25}
       There exists a curve\/ $c:((-\ep,\ep),0)\to(\R^2,\zv)$ such that\/
       $c'(0)$ is parallel to\/ $\eta(\zv)$,
       $\hat{c}'_0=0, \hat{c}''_0\ne0$ and there exists\/ $\ell$
       satisfying\/ $\hat{c}'''_0=\ell\hat{c}''_0$ and\/
       $a:=
       \det(\hat{\gamma}',
        \hat{c}'',
        3\hat{c}^{(5)}-10\ell\hat{c}^{(4)}
       )(0)\ne0$, where\/
       $\hat{c}=f\circ c$ and\/ $\hat{\gamma}=f\circ \gamma$.
  \item\label{it:contact}
       $\psi(0)=\psi'(0)=\cdots=\psi^{(k-1)}(0)=0$ and\/
       $b:=\psi^{(k)}(0)\ne0$, where\/
       $\psi$ is the function defined by\/ {\rm (\ref{eq:ccr})}.
  \item\label{it:sign}
       If\/ $k$ is even, sign\/ $\pm$ of the cuspidal\/ $S_{k-1}^\pm$
       singularity
       coincides with the sign of the product \/
       $ab=\det(\hat{\gamma}',
        \hat{c}'',
        3\hat{c}^{(5)}-10\ell\hat{c}^{(4)}
       )(0)\cdot\psi^{(k)}(0)$. Here, we choose\/ $\eta$ and\/ $t$
       so that\/
       $c'(0)$ points the same direction
       as the null vector\/ $\eta(0)$ and that\/
       $(\gamma',\eta)(0)$ is positively oriented.
 \end{enumerate}
\end{thm}
To prove Theorem \ref{thm:csk}, we show at first
 the following lemma.
\begin{lem}
 \label{independ}
 The conditions in Theorem\/ {\rm \ref{thm:csk}} do not depend either
 choice of coordinates on the source, the parameter of\/ $c$,
 the parameter of\/ $\gamma$, the choice of representative\/ $\nu$,
 the choice of\/ $\eta$
 or on the 
 choice of coordinates on the target.
\end{lem}
It is easy to check that the condition \eqref{it:ndg} does not depend on
all the choices by Lemma \ref{lem:25}.
Since linear independence is not changed
by a diffeomorphism, the condition \eqref{it:25} does not depend
on all the choices.
We shall prove that either of the
conditions \eqref{it:contact} and \eqref{it:sign}
does not depend on all the choices.
\begin{proof}[Proof for the condition\/ \eqref{it:contact}.]
 Note that the condition \eqref{it:contact} is not changed
 on the non-zero functional multiple of $\psi$ on $S(f)$.
 Thus it does not depend on the choices of
 $\nu$, $\eta$ and the parameter of $\gamma$.
 Hence it is sufficient to prove that the condition \eqref{it:contact}
 does not depend on the choice of the coordinates on the target.

 Let $\Phi:(\R^3,\zv)\to(\R^3,\zv)$ be a diffeomorphism-germ
 and $d\Phi$ its derivative.
 The map $d\Phi$ can be considered as a $GL(3,\R)$-valued
 function $q\mapsto d\Phi_{q}$.
 Since $Au\times Av=(\det A)\,{}^t{A^{-1}}u\times v$ for any
 vectors $u$ and $v$ in the 3-space and any non-singular
 matrix $A$, we can take
 $\tilde\nu={}^t(d\Phi)^{-1}\nu$
 as a normal vector field of $\tilde{f}=\Phi\circ f$.
 So, we shall prove
 $$\tilde\psi(t)=
 \det\big((\tilde{f}\circ\gamma)'(t),\ 
 \tilde{\nu}\circ\gamma(t),\ 
 d\tilde\nu_{\gamma(t)}(\eta(t))\big)
 $$
 is a non-zero functional multiple of $\psi(t)$.
 
 Since the condition does not depend on the choices of
 coordinates on the source, choice of $\eta$ 
 and choice of $\nu$,
 we may assume that $S(f)=\{v=0\}$, $\eta=\partial/\partial v$
 on $\gamma(t)$,
 $\nu$ is the unit normal vector and
 $f(u,0)$ is the arc-length parameter.
 Under this assumption,
 $f_u,\nu,\nu_v$ are orthogonal each other,
 since we see $\inner{\nu}{\nu_v}=0$ from $\inner{\nu}{\nu}=1$
 and
 $\inner{f_u}{\nu_v}=-\inner{f_{uv}}{\nu}=\inner{f_v}{\nu_u}=0$
 on $S(f)$ from $f_v=\zv$ on $S(f)$.
 Hence $\nu\times\nu_v$ is parallel to $f_u$.
 Thus we have
 $\psi=\det(f_u,\nu,\nu_v)=\inner{f_u}{\nu\times\nu_v}$,
 and it holds that $\nu\times\nu_v=\psi f_u$.
 Then, it follows that
 $$
 \begin{array}{rcl}
 \tilde\psi(t)&=&
 \det
 \Big(
 d\Phi_{f(\gamma(t))}f_u(\gamma(t)),\ 
 {}^t(d\Phi_{f(\gamma(t))})^{-1}\nu(\gamma(t)),\ 
 \big({}^t(d\Phi_{f(\gamma(t))})^{-1}
 \nu(\gamma(t))\big)_v\Big)\\[3mm]
 &=&
 \det\big(
 d\Phi_{f(\gamma(t))}f_u(\gamma(t)),\ 
 {}^t(d\Phi_{f(\gamma(t))})^{-1}\nu(\gamma(t)),\ 
 {}^t(d\Phi_{f(\gamma(t))})^{-1}\nu_v(\gamma(t))
 \big).
\end{array}
 $$
 Here note that 
 $({}^t(d\Phi_{f(\gamma(t))})^{-1})_v=0$ on $S(f)$ 
 because $f_v(\gamma(t))=0$ and
 $\det(d\Phi_{f(\gamma(t))})\neq 0$.
 Omitting $(t)$, $\gamma(t)$ and $f(\gamma(t))$,
 we can modify
 \begin{equation*}
    \begin{array}[tb]{rcl}
     \tilde\psi&=&\det
      \big(
      d\Phi f_u,\ {}^t(d\Phi)^{-1}\nu,\ {}^t(d\Phi)^{-1}\nu_v\big)
      =\big\langle{d\Phi f_u},\ 
       {{}^t(d\Phi)^{-1}\nu\times{}^t(d\Phi)^{-1}\nu_v}\big\rangle\\
      &=&\big\langle{d\Phi f_u},\ 
         {\det({}^t(d\Phi)^{-1})d\Phi (\nu\times\nu_v)}\big\rangle
      =\det({}^t(d\Phi)^{-1})\inner{d\Phi f_u}{d\Phi (\psi f_u)}\\
     &=&\det((d\Phi)^{-1})\inner{d\Phi f_u}{d\Phi f_u}\psi.\\
 \end{array}\end{equation*}
 Since $\det((d\Phi)^{-1})\inner{d\Phi f_u}{d\Phi f_u}$ is a function
 which never vanishes on $S(f)$,
 the condition \eqref{it:contact} does not depend on the choice of the
 coordinate system on the target.
\end{proof}

\begin{proof}[Proof for the condition\/ \eqref{it:sign}.]
 When the direction of the representative $\nu$
 of $[\nu]$ is changed to opposite direction,
 the signs of both $a$ and
 $b$ are not changed.
 When the parameter of $\gamma$ reverses,
 the sign of $a$ is unchanged,
 and if $k$ is even then the sign of $b$ is unchanged
 because of the positivity of the basis $(\gamma',\eta)$.
 If the orientation of the target is changed,
 then signs of both $a$ and $b$ are changed.
 Hence in all the cases, $\sign(ab)$ is not changed.
\end{proof}
\begin{proof}[Proof of Theorem\/ {\rm \ref{thm:csk}}]
 Assume that a map-germ $f$ satisfies the conditions in
 Theorem \ref{thm:csk}.
 By the same argument as in the proof of Theorem \ref{thm:cmm},
 we may assume that
 $$f(u,v)=\big(u,vg(u,v),vh(u,v)\big).$$
 Consider the new coordinate system $(u,\tilde v)$ satisfying
 $S(f)=\{\tilde v=0\}$ and rewrite $\tilde{v}$ by $v$.
 Then,we get $g=h=0$ on $v=0$.
 Thus, there exist functions $\tilde g(u,v)$ and
 $\tilde h(u,v)$ such that
 $f(u,v)=\big(u,v^2\tilde g(u,v),v^2\tilde h(u,v)\big).$
 By the same argument as in the proof of Theorem \ref{thm:cmm}
 again, we may assume that
 $$f(u,v)=(u,v^2,v^3\overline{h}(u,v^2)).$$
 The normal vector $[\nu]$ of $f$  and $d\nu(\eta)$ are
 given by
 $$
 \nu=\left(*, 
 -3v\bar{h}(u,v^2)-2v^3\frac{\partial \bar{h}}{\partial
 v}(u,v^2),2\right)\text{\ and\ }
 d\nu(\eta)=d\nu(\partial/\partial v)
 =\left(*,-3\bar{h}(u,v^2)+v*,0\right).
 $$
 Thus the function $\psi$ of this map is $6\overline{h}(t,0)$.
 Thus the condition (c) is written as
 $$\overline{h}=
 (\partial/\partial u)\overline{h}=
 \ldots=
 (\partial^{k-1}/\partial u^{k-1})\overline{h}=0 \text{\ and\ }
  (\partial^k/\partial u^k)\overline{h}\ne0
 \text{\ at\ }\zv.$$
 Consider a curve $(o(v),v)$,
 where $o(v)$ is a Landau notation again.
 Then this curve is tangent to
 $\eta$ at $\zv$ and all curves passing through
 $\zv$ tangent
 to $\eta=\partial/\partial v$ at
 $\zv$ are written by this form.
 Since $(\partial/\partial u)\overline{h}(\zv)=0$,
 $v^3\overline{h}(o(v),v)$ has no terms of $v^4$.
 Thus the condition (b) is equivalent to
 $(\partial^5/\partial
 v^5)v^3\bar{h}(o(v),v)(0)\ne0$.
 Hence, the coefficient of $v^2$
 in $\overline{h}(u,v^2)$ is not zero.
 Thus it follows that there exist functions $h,\ \overline{g}$
 and non-zero real numbers $\alpha,\beta$ such that
 $$
 \overline{h}(u,v^2)
 =
 \alpha v^2+\beta u^kh(u)+v^2\alpha\overline{g}(u,v^2),
 \quad h(0)=1.
 $$
 By the coordinate system change
 \begin{equation}
  \label{eq:change1}
   \begin{array}[tb]{l}
    U=u\sqrt[k]{h(u)}\\
    V=\sqrt{|\alpha|}v\sqrt{1+\overline{g}(u,v^2)},
   \end{array}
 \end{equation}
 $\overline{h}$ becomes
 $\sign(\alpha)V^2+\beta U^k$.
 One can easily see that
 the inverse map of \eqref{eq:change1} is given by
 $$
 \begin{array}[tb]{l}
  u=UH(U)\\
  v=VG(U,V^2),
 \end{array}
 $$
 using functions $G,H$ whose constant terms are not zero.
 Hence $f$ is ${\cal A}$-equivalent to
 $$
 f(U,V)=
 \big(UH(U),V^2G(U,V^2)^2,V^3G(U,V^2)^3(\sign(\alpha)V^2+\beta U^k)\big).
 $$
 Now we consider a map-germ
 $\big(u,v^2,v^3(\sign(\alpha)v^2+\beta u^k)\big)$ and
 a diffeomorphism $$\Psi(X,Y,Z)=\big(XH(X),YG(X,Y)^2,ZG(X,Y)^3\big).$$
 Then it follows that $f$ is ${\cal A}$-equivalent to
  $\big(u,v^2,v^3(\sign(\alpha)v^2+\beta u^k)\big)$.

 Here, we have $ab=6(6!k!)\sign(\alpha)\beta$.
 By a suitable scale change,
 if $k$ is odd or $k$ is even and $\sign(\alpha)\beta>0$, then
 $f$ is ${\cal A}$-equivalent to $(u,v^2,v^3(v^2+u^k))$.
 If $k$ is even and $\sign(\alpha)\beta<0$, then 
 $f$ is ${\cal A}$-equivalent to $(u,v^2,v^3(v^2-u^k))$.
\end{proof}
\section{Applications}
\label{sec:app}
In this section,
we give two applications of our criteria.
\medskip

\begin{figure}[htbp]
\centering
   \includegraphics[width=.7\linewidth]{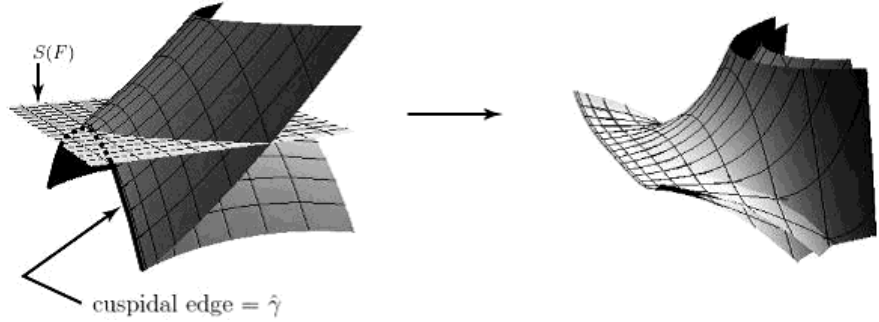}
 \caption{Left : $2$-point contact of $\hat\gamma$ and $S(F)$\quad
  Right : cuspidal $S_1^-$ singularity of $F\circ f$}
 \label{fig:s2}
 \end{figure}

Let $\s:((-\ep,\ep),0)\to(\R^3,\zv)$ be a space curve such that
its curvature
never vanishes, with the
arclength parameter.
Let $\e,\n,\b$ be its Frenet frame and
$\kappa$, $\tau$ its
curvature and torsion respectively.
A map $(t,u)\mapsto \s(t)+u\e(t)$
is called the {\it tangent developable surface} of $\s$.
In \cite{Mond3}, Mond proved the following theorem.
\begin{thm}[Mond \cite{Mond3}]
 \label{thm:tang}
 The germ of the tangent developable surface
 of\/ $\s$
 at\/ $(0,0)$ is\/ ${\cal A}$-equivalent
 to the cuspidal\/ $S_{1}^+$ singularity
 if\/ $\tau=\tau'=0$ and\/ $\tau{''}\ne0$ at $0$.
\end{thm}
\begin{rem}
 Mond also classified the case that\/
 $\tau=\tau'=\cdots=\tau^{(k-1)}=0$, $\tau^{(k)}\ne0$ for\/
 $k=3$ and\/ $4$. By Ishikawa's theorem\/ {\rm  \cite{goo}},
 the developable surfaces do not have any cuspidal\/ $S_{k-1}$
 singularities for\/ $k>2$.
\end{rem}
We shall prove Theorem \ref{thm:tang}
using our criteria as an application.
\begin{proof}
 Let $\s$ be a space curve and
 $f(t,u)=\s(t)+u\e(t)$ the
 tangent developable surface of $\s$.
 Then $S(f)=\{u=0\}$ and
 $\eta=-\partial/\partial t+\partial/\partial u$.
 Since $\lambda=\det(\e+u\kappa \n,\e,\b)=-\kappa u$, we
 see that $d\lambda\ne0$ and the singularities are non-degenerate.
 Let us consider a curve
 $$
 c:t\mapsto\left(-t,-\frac{\s(-t)\cdot\e(0)}
 {\e(-t)\cdot\e(0)}\right),
 $$
  in the $(t,u)$-space and put $\hat{c}=f\circ c$.
 Then,
 we see that $c$ satisfies the
 condition \eqref{it:25} of
 Theorem \ref{thm:csk}.
 In fact, by a direct calculation
 we have
 $$
 \begin{array}[tb]{l}
  \hat{c}'(0)=\zv,\quad
   \hat{c}''(0)=-\kappa(0)\n(0),\quad
   \hat{c}'''(0)=2\kappa'(0)\n(0),\\
  \hat{c}^{(4)}(0)=*\e(0)+*\n(0),\quad\text{and}\quad
   \hat{c}^{(5)}(0)=*\e(0)+*\n(0)+4\kappa(0)\tau''(0)\b(0).
 \end{array}
 $$
 Hence
 $a=-12\det(\e,\kappa\n,\kappa\tau''\b)(0)\ne0$ holds.
 Moreover, since we can take $\nu=\b$, we have
 $d\nu(\eta)=-(\partial \b/\partial t)=\tau\n$.
 So,
 $$
 \psi(t)=\det\big(\hat{\gamma}'(t),\ \nu(\gamma(t)),\
 \eta(t)\nu(\gamma(t))\big)
 =
 \det(\e(t),\b(t),\tau \n(t))=-\tau(t)
 $$
 and hence $b=-\tau''(0)$ by the assumption.
 Since $ab>0$, $f$ at $(0,0)$ is ${\cal A}$-equivalent to 
 the cuspidal $S_{1}^+$ singularity by Theorem \ref{thm:csk}.
\end{proof}
Now we consider another property of the cuspidal $S_k$ singularity.
The following observation about the cuspidal cross
cap was given by Arnol'd \cite[p.120 Example 3]{arnld-sing-caus}.
Let $f:(\R^2,\zv)\to(\R^3,\zv)$ be a cuspidal edge and
$F:\R^3\to\R^3$ a ``generic'' fold.
Then the map-germ $F\circ f$ at $\zv$ is a cuspidal cross cap,
where the {\it cuspidal edge\/} is a map-germ defined by
$(u,v)\mapsto(u,v^2,v^3)$ at the origin and
the {\it fold\/} is a map-germ defined by
$(x,y,z)\mapsto (x,y,z^2)$ at the origin.
Here, we generalize this observation and clarify
the meaning of {\it genericness}.
 \begin{thm}
  \label{thm:app2}
  Let\/ $f:(\R^2,\zv)\to(\R^3,\zv)$ be a map-germ\/ ${\cal A}$-equivalent to
  the cuspidal edge
  and\/ $F:(\R^3,\zv)\to(\R^3,\zv)$ be a map-germ\/ ${\cal A}$-equivalent to
  the fold.
  Assume that the following three conditions:
  \begin{itemize}
   \item[{\rm (A)}] the limiting tangent plane\/ $LT$ of\/
                $f$ at\/ $\zv$ does not contain
                the kernel\/ $\ker dF_{\zv}$,
   \item[{\rm (B)}] $LT$ is transverse to\/ $S(F)$ at\/ $\zv$,
   \item[{\rm (C)}] the singular curve\/ $\hat\gamma=f(S(f))$ has
    a\/ $k$-point contact with\/ $S(F)$ at\/ $\zv$.
  \end{itemize}
  Then the
  composition\/ $F\circ f$ at\/ $\zv$ is\/
  ${\cal A}$-equivalent to the cuspidal\/ $S_{k-1}^\pm$ singularity.

  In the case of\/ $k$ is even, the sign\/ $\pm$ is determined
  by the following rule.
  Since\/ $k$ is the order of contact
  between\/ $\hat\gamma$ and\/ $S(F)$,
  it holds that\/ $\hat\gamma$ is locally located on the half
  space bounded by\/ $S(F)\subset\R^3$.
  If\/ $\Im(f)$ is also locally located only
  on the same side as\/ $\hat\gamma$,
  then\/ $F\circ f$ is\/ ${\cal A}$-equivalent 
  to the cuspidal\/ $S_{k-1}^+$
  singularity.
  Otherwise, $F\circ f$ is\/ ${\cal A}$-equivalent 
  to the cuspidal\/ $S_{k-1}^-$
  singularity\/ $($See Figure\/ {\rm \ref{fig:s2}}$)$.
 \end{thm}
 Since the 
 condition (A) means that the normal vector $\nu(\zv)$ of $f$ is not
 perpendicular to $\ker dF_{\zv}$
 and the condition (B) means that the normal vector $\nu(\zv)$ of $f$
 does not perpendicular to
 the tangent plane of $S(F)$,
 the conditions (A) and (B) are generic conditions.
 It should be remarked that
 folding maps for smooth surfaces
 are considered in \cite{bw,it}.
\begin{proof}
 Let $f$ and $F$ be map-germs satisfying
 the conditions in Theorem \ref{thm:app2}.
 In the following diagram:
 $$
 \begin{CD}
  (\R^2;(u,v),\zv)@>{f}>>(\R^3;(x,y,z),\zv)@>{F}>>(\R^3;(X,Y,Z),\zv),
 \end{CD}
 $$
 the conditions and assertions of
 the theorem do not depend on the choice of
 coordinate systems on each space.
 So, we take the coordinate systems $(x,y,z)$ and $(X,Y,Z)$
 which satisfy
 $F(x,y,z)=(x,y,z^2)$.
 Moreover, we can take the coordinate system $(u,v)$ so
 that $S(f)=\{v=0\}$ and $\eta=\partial/\partial v$ there.
 Denote $f(u,v)=(f_1(u,v),f_2(u,v),f_3(u,v))$.
 Then by the transversality condition, 
 either $(\partial/\partial u) f_1(0,0)\ne0$ or
 $(\partial/\partial u) f_2(0,0)\ne0$ holds.
 By the implicit function theorem,
 we may assume $f(u,v)=(u,f_2(u,v),f_3(u,v))$.
 Then by the conditions $S(f)=\{v=0\}$ and $\eta=\partial/\partial v$,
 $f$ has the 
 following form: $f(u,v)=(u,a_2(u)+v^2b_2(u,v),a_3(u)+v^2b_3(u,v))$.
 By the condition (A) that the
 limiting tangent plane does not contain
 the $Z$-axis,
 it holds that $b_2(0,0)\ne0$.
 By the coordinate change
 $\tilde{u}=u,\ \tilde{v}=v\sqrt{b_2(u,v)}$,
 we may assume $f(u,v)=(u,a_2(u)+v^2,a_3(u)+v^2b_3(u,v))$.
 Since $\alpha(x,y,z)=(x,y-a_2(x),z)$ and
 $\beta(X,Y,Z)=(X,Y-a_2(X),Z)$ are both diffeomorphism,
 considering $\alpha\circ f$ and $\beta\circ F\circ\alpha^{-1}$,
 and rewriting them by $f$ and $F$, we may assume that
 $f(u,v)=(u,v^2,a_3(u)+v^2b_3(u,v))$ and $F(x,y,z)=(x,y,z^2)$
 again.
 So,
 $$
 F\circ f(u,v)=
 \big(u,v^2,a_3(u)^2+2v^2a_3(u)b_3(u,v)+v^4b_3(u,v)^2\big).
 $$
 Then $\partial (F\circ f)/\partial u=(1,0,2a_3a_3'+2v^2\,*)$
 and $\partial (F\circ f)/\partial v=2v(0,1,2a_3b_3+va_3(b_3)_v+v^2\,*)$,
 where $*$ means a function, $a_3'=d a_3/du$
 and $(b_3)_v=\partial b_3/\partial v$.
 We can take $\nu=(2a_3a_3'+2v^2\,*,2a_3b_3+va_3(b_3)_v+v^2\,*,-1)$
 as a normal vector.
 Then we see that the signed area density function $\lambda$
 is a non-zero functional multiple of $v$.
 Thus the non-degeneracy of all singularities
 of $F\circ f$ follows.
 Since $f|_{S(f)}=(u,0,a_3(u))$,
 the condition (C) implies
 $a_3'=\cdots=a_3^{(k-1)}=0$ and $a_3^{(k)}\ne0$ at $0$.
 
 Since $f$ is ${\cal A}$-equivalent to the cuspidal edge at the origin,
 $(b_3)_v(0,0)\ne0$ holds.
 The function $\psi$ of
 $F\circ f$ defined by (\ref{eq:ccr}) is given by
 $3a_3(t)(b_3)_v(t,0)(1+4a_3(t)^{2}a_3'(t)^{2})$
 because of
 $d\nu(\eta)(t,0)=\nu_v(t,0)
 =\big(0,3a_3(t)(b_3)_v(t,0),0\big)$.

 If $k=1$, then we have the conclusion
 by Corollary 1.5 of \cite{FSUY}.
 If $k\geq2$, we have $b_3(0,0)\ne0$ by the
 transversality condition (B).
 Now, we consider a curve $c(t)=(0,t)$
 and $\hat{c}(t)=F\circ f(c(t))=\big(0,t^2,t^4b_3(0,t)^2\big)$.
 Since $b_3(0,0)\ne0$ and $(b_3)_v(0,0)\ne0$,
 we see that the conditions \eqref{it:25} of Theorem \ref{thm:csk}
 are satisfied.
 If $k$ is odd, by Theorem \ref{thm:csk} $F\circ f$ is ${\cal
 A}$-equivalent
 to the cuspidal $S_{k-1}$ singularity.
 If $k$ is even, $F\circ f$ is equivalent to
 the cuspidal $S_{k-1}^+$ singularity (resp. the
 cuspidal $S_{k-1}^-$ singularity) 
 if and only if $a_3^{(k)}(0)b_3(0,0)>0$
 (resp. $a_3^{(k)}(0)b_3(0,0)<0$).
 Since $S(F)=\{(x,y,z)\,|\,z=0\}$
 and $f(u,v)=(u,v^2,a_3(u)+v^2b_3(u,v))$,
 one can easily see that the condition $a_3^{(k)}(0)b_3(0,0)>0$
 is equivalent to the condition 
 that $\Im(f)$ is locally located on the half space
 bounded by the $xy$-plane such that 
 $f(u,0)$ lies in. 
 This completes the proof.
\end{proof}
Mond's criteria \cite[Theorem 4.1.1]{Mond2}
is useful for normalized germs.
But in general, like for the examples of
this section, our criteria seem more useful.

\def\cprime{$'$}

\begin{flushright}
  \begin{tabular}[h]{l}
    Department of Mathematics,\\
    Faculty of Education,\\
    Gifu University,\\
    Yanagido 1-1, Gifu, 501-1193, Japan.\\
    {\tt ksajiO$\!\!\!$agifu-u.ac.jp}
  \end{tabular}
\end{flushright}
\end{document}